\documentclass[12pt]{article}
\usepackage{amsfonts,amsmath,amsxtra}
\usepackage{latexsym}
\usepackage{amssymb}
\usepackage{color}
\usepackage[matrix,arrow,curve]{xy}

\makeatletter
\def\@seccntformat#1{\csname
the#1\endcsname\enspace} \makeatother

\def\hybrid{\topmargin 0pt      \oddsidemargin 0pt
        \headheight 0pt \headsep 0pt
        \textwidth 16.5cm
        \textheight 23cm
        \marginparwidth 0.0in
        \parskip 5pt plus 1pt   \jot = 1.5ex}
\catcode`\@=11
\def\marginnote#1{}
\newcount\hour
\newcount\minute
\newtoks\amorpm
\hour=\time\divide\hour by60 \minute=\time{\multiply\hour by60
\global\advance\minute by-\hour}
\edef\standardtime{{\ifnum\hour<12 \global\amorpm={am}%
        \else\global\amorpm={pm}\advance\hour by-12 \fi
        \ifnum\hour=0 \hour=12 \fi
      \number\hour:\ifnum\minute<10 0\fi\number\minute\the\amorpm}}
\edef\militarytime{\number\hour:\ifnum\minute<10 0\fi\number\minute}

\def\draftlabel#1{{\@bsphack\if@filesw {\let\thepage\relax
   \xdef\@gtempa{\write\@auxout{\string
      \newlabel{#1}{{\@currentlabel}{\thepage}}}}}\@gtempa
   \if@nobreak \ifvmode\nobreak\fi\fi\fi\@esphack}
        \gdef\@eqnlabel{#1}}
\def\@eqnlabel{}
\def\@vacuum{}
\def\draftmarginnote#1{\marginpar{\raggedright\scriptsize\tt#1}}

\def\draft{\oddsidemargin -0.1truein
        \def\@oddfoot{\sl HeckeBaxter.tex \hfil
        \rm\thepage\hfil\sl\today\quad\militarytime}
        \let\@evenfoot\@oddfoot \overfullrule 3pt
        \let\label=\draftlabel
        \let\marginnote=\draftmarginnote
\def\@eqnnum{{\rm (\theequation)}
\rlap{\kern\marginparsep\tt\@eqnlabel}%
\global\let\@eqnlabel\@vacuum}  }

\newfont{\Bbbb}{msbm7 scaled 1\@ptsize00}
\newcommand{\zs}{\raise-1pt\hbox{$\mbox{\Bbbb Z}$}}

scaled 1\@ptsize00

\def\numberbysection{\@addtoreset{equation}{section}
        \def\theequation{\thesection.\arabic{equation}}}
\numberbysection

\renewcommand{\theequation}{\thesection.\arabic{equation}}

\def\titlepage{\@restonecolfalse\if@twocolumn\@restonecoltrue\onecolumn
     \else \newpage \fi \thispagestyle{empty}\c@page\z@
\def\thefootnote{\fnsymbol{footnote}} }
\def\endtitlepage{\if@restonecol\twocolumn \else  \fi
        \def\thefootnote{\arabic{footnote}}
        \setcounter{footnote}{0}}  
\relax \hybrid
\makeatletter
\newdimen\normalarrayskip            
\newdimen\minarrayskip               
\normalarrayskip\baselineskip \minarrayskip\jot
\newif\ifold             \oldtrue            \def\new{\oldfalse}
\def\arraymode{\ifold\relax\else\displaystyle\fi}
\def\eqnumphantom{\phantom{(\theequation)}} 
\def\@arrayskip{\ifold\baselineskip\z@\lineskip\z@
     \else
     \baselineskip\minarrayskip\lineskip1\baselineskip\fi}
\def\@arrayclassz{\ifcase \@lastchclass \@acolampacol \or
\@ampacol \or \or \or \@addamp \or
   \@acolampacol \or \@firstampfalse \@acol \fi
\edef\@preamble{\@preamble
  \ifcase \@chnum
     \hfil$\relax\arraymode\@sharp$\hfil
     \or $\relax\arraymode\@sharp$\hfil
     \or \hfil$\relax\arraymode\@sharp$\fi}}
\def\@array[#1]#2{\setbox\@arstrutbox=\hbox{\vrule
     height\arraystretch \ht\strutbox
     depth\arraystretch \dp\strutbox
width\z@}\@mkpream{#2}\edef\@preamble{\halign \noexpand\@halignto
\bgroup \tabskip\z@ \@arstrut \@preamble \tabskip\z@ \cr}%
\let\@startpbox\@@startpbox \let\@endpbox\@@endpbox
  \if #1t\vtop \else \if#1b\vbox \else \vcenter \fi\fi
  \bgroup \let\par\relax
  \let\@sharp##\let\protect\relax
  \@arrayskip\@preamble}
\def\eqnarray{\stepcounter{equation}%
              \let\@currentlabel=\theequation
              \global\@eqnswtrue
              \global\@eqcnt\z@
              \tabskip\@centering              
              \let\\=\@eqncr
              $$%
            \halign to \displaywidth  \bgroup
             \eqnumphantom \@eqnsel
      \hskip\@centering                               
    $\displaystyle  \tabskip\z@ {##}$%
    &\global\@eqcnt\@ne \hskip 2\arraycolsep
         $ \displaystyle  \arraymode{##}$\hfil
    &\global\@eqcnt\tw@ \hskip 2\arraycolsep
         $\displaystyle\tabskip\z@{##}$\hfil
         \tabskip\@centering
    &{##}\tabskip\z@\cr}
\makeatother

\def\IC{\mathbb{C}}

\def\IQ{\mathbb{Q}}
\def\IR{\mathbb{R}}
\def\IZ{\mathbb{Z}}
\def\CA {\mathcal{A}}

\def\CD {\mathcal{D}}

\def\CH {\mathcal{H}}

\def\CS {\mathcal{S}}

\def\CU {\mathcal{U}}
\def\CV {\mathcal{V}}
\def\CW {\mathcal{W}}

\def\fh{{\frak h}}

\def\fk{{\frak k}}

\def\fu{{\frak u}}

\def\g {{\gamma}}

\def\pr {\partial}

\def\diag{{\rm diag}}

\def\End{{\rm End}}
\def\Id{{\rm Id}}

\def\Lie{{\rm Lie}}
\def\Mat{{\rm Mat}}

\def\Span{{\mathop{\rm span}}}

\def\frak{\mathfrak}
\def\ov {{\overline}}

\def\Tr{{\rm Tr}\,}

\def\gl{\mathfrak{gl}}
\def\hgl{\mathfrak{hgl}}

\def\ssp{\mathfrak{sp}}
\def\hsp{\mathfrak{hsp}}

\def\<{\langle}
\def\>{\rangle}
\def\ov{\overline}

\makeatletter
\DeclareRobustCommand{\loplus}{\mathbin{\mathpalette\dog@lsemi{+}}}
\DeclareRobustCommand{\roplus}{\mathbin{\mathpalette\dog@rsemi{+}}}

\newtheorem{te}{Theorem}[section]
\newtheorem{de}{Definition}[section]

\newtheorem{lem}{Lemma}[section]

\newcommand{\proof}{\noindent {\it Proof}. }
\newcommand\bqa{\begin{eqnarray}}
\newcommand\eqa{\end{eqnarray}}
\def\be{\begin{eqnarray}\new\begin{array}{cc}}
\def\ee{\end{array}\end{eqnarray}}
\def\beq{\begin{equation}}
\def\eeq{\end{equation}}
\def\bse{\begin{subequations}}                
\def\ese{\end{subequations}}
\def\bp{\begin{pmatrix}}
\def\ep{\end{pmatrix}}

\def\i{\imath}
\newcounter{pac}[section]

\newcounter{pacc}[subsection]

 0

\begin{document}

\title{\bf The Hecke-Baxter operators via\\  Heisenberg group extensions}
\author{A.A. Gerasimov, D.R. Lebedev and S.V. Oblezin}
\date{\today}
\maketitle

\renewcommand{\abstractname}{}

\begin{abstract}
\noindent {\bf Abstract}. The $GL_{\ell+1}(\IR)$ Hecke-Baxter
operator was introduced as an element of the $O_{\ell+1}$-spherical
Hecke algebra associated with the Gelfand pair $O_{\ell+1}\subset
GL_{\ell+1}(\IR)$. It was specified by the property to act on an
$O_{\ell+1}$-fixed vector in a $GL_{\ell+1}(\IR)$-principal series
representation via multiplication by the local Archimedean
$L$-factor canonically attached to the representation. In this note
we propose another way to define the Hecke-Baxter operator,
identifying it with a generalized Whittaker function for an
extension of the Lie group $GL_{\ell+1}(\IR)\times GL_{\ell+1}(\IR)$ by a Heisenberg Lie
group.  We also show how this  Whittaker function can be lifted
  to a matrix element of an extension of the Lie group
  $Sp_{2\ell+2}(\IR)\times Sp_{2\ell+2}(\IR)$ by a Heisenberg Lie
group.

\end{abstract}
\vspace{5mm}

\section{Introduction}

Introducing  the $GL_{\ell+1}(\IR)$ Hecke-Baxter operator  in
\cite{GLO08} is a  result of combining two different lines of
research. Along the one direction, to each admissible representation
of the split real Lie  group $GL_{\ell+1}(\IR)$ one might
canonically attach a local Archimedean $L$-factor, a function in one
complex variable capturing data of the considered representation
(see e.g. \cite{ILP} for pedagogical introduction). Consideration of
this particular function is motivated by deep connections with the
theory of zeta-functions associated with global number fields. Also
Archimedean  $L$-factors attached to representations of non-compact
real Lie groups play an essential role in formulation of the
Langlands correspondence for such groups. However the standard
construction of $L$-factor for a given representation followed a
very round-about strategy pioneered in \cite{JL}.

Meanwhile, along quite different line of research R. Baxter
introduced a special type of operators acting in the space of states
of lattice quantum integrable systems \cite{Ba}. This type of
operators turns out to be a particular instance  of a universal
construction instrumental in finding  exact solutions of various
quantum integrable systems. For example, in \cite{PG} the Baxter
operator is proposed for the closed Toda chain. In spite of large
amount of work, the underlying reason for considering these
particular operators  is not manifested  (which is quite usual case
with {\it ad hoc} constructions in the theory of quantum integrable
systems).

Fortunately, these two lines of development intersect in the theory
of $GL_{\ell+1}(\IR)$-Toda chains, one of the simplest but still
interesting exactly solvable quantum mechanical systems having
numerous ties with various branches of mathematics and mathematical
physics. Solutions of $GL_{\ell+1}(\IR)$-Toda chains (common
eigenfunctions of the mutually commuting quantum Hamiltonians)
 may be expressed in terms of a particular set of matrix elements in principal
series spherical representations of $GL_{\ell+1}(\IR)$. These matrix
elements appear to be generalizations of a particular class of the
Whittaker functions \cite{Ja}. This  quantum integrable system also
allows explicit construction of the Baxter operator as an integral
operator. In this regard in \cite{GLO08} (see also \cite{G}), a
particular element of the spherical Hecke algebra
$\CH(GL_{\ell+1}(\IR),O_{\ell+1})$ associated with the Gelfand pair
$O_{\ell+1}\subset GL_{\ell+1}(\IR)$ is introduced. It is
demonstrated that its action on  the $GL_{\ell+1}(\IR)$-Whittaker
functions descends to the action of the Baxter integral operator for
the $GL_{\ell+1}$-Toda chain. Moreover the Whittaker function is an
eigenfunction of this integral operator with the eigenvalue equal to
the local Archimedean $L$-factor attached to the considered
principal series $GL_{\ell+1}(\IR)$-representation. This result is
not restricted to the Whittaker functions, since the action of this
particular element of the spherical Hecke algebra on the spherical
vector in principal series representation is reduced to
multiplication by the corresponding $L$-factor. We coined the term
Hecke-Baxter operator for this element of spherical Hecke algebra
$\CH(GL_{\ell+1}(\IR),O_{\ell+1})$.

Although general framework for constructing the Baxter operators in
representation theory started to emerge, the question of why we
should consider this particular element of the spherical Hecke
algebra is yet to be clarified. Its relation with local Archimedian
$L$-factor is intriguing but does not substantiate the appearance of
this particular element of spherical Hecke algebra as the
Hecke-Baxter operator (the theory of $L$-factors itself is a kind of
mystery and requires further clarifications).

In this note we demonstrate that the $O_{\ell+1}$-biinvariant
spherical function on $GL_{\ell+1}$ entering the description of the
Hecke-Baxter operator may be identified with a (generalized)
Whittaker function on a non-reductive Lie group obtained as an
extension $\CH GL_{\ell+1}(\IR)$ of $GL_{\ell+1}(\IR)\times GL_{\ell+1}(\IR)$ by a
Heisenberg Lie group. Precisely, generalizing  the definition of Whittaker function to the
case of non-reductive Lie groups, we construct the Whittaker
function as a matrix element of an appropriate irreducible
representation of the Lie group extension $\CH GL_{\ell+1}(\IR)$ of
$GL_{\ell+1}(\IR)\times GL_{\ell+1}(\IR)$ by a Heisenberg Lie group. The corresponding
Heisenberg Lie algebra $\CH^{(\ell+1)}$ is the standard Heisenberg
algebra associated with the vector space $\Mat_{\ell+1}(\IR)$ of
real $(\ell+1)\times (\ell+1)$ matrices. The right and left actions of
two copies of $GL_{\ell+1}(\IR)$ on $\Mat_{\ell+1}(\IR)$ provides the
extension $\CH GL_{\ell+1}(\IR)$ of $GL_{\ell+1}(\IR)\times GL_{\ell+1}(\IR)$ by the
Heisenberg Lie group. The resulting extension may be called the
matrix hyperbolic/inverted harmonic oscillator group, since in the
case $\ell=0$ this group reduces to the hyperbolic version of the
harmonic oscillator group (see e.g. \cite{K},\cite{Ho}).

Our construction of Whittaker functions for $\CH GL_{\ell+1}$ has
  some obvious hidden symmetries. This  suggests that the
    considered Whittaker
 functions might be naturally lifted to  the Whittaker functions
   given by the matrix elements
 of larger  Lie groups. Indeed at the end of this note we construct Whittaker
 functions for a Heisenberg group extension of $Sp_{2\ell+2}\times
 Sp_{2\ell+2}$ and demonstrate that a special instance of these
 functions reproduce the Whittaker functions for $\CH
 GL_{\ell+1}(\IR)$ constructed earlier.
Thus representation theory of symplectic Lie groups and their
Heisenberg extensions looks like a natural habitat for the
$GL_{\ell+1}(\IR)$ Hecke Baxter operator.

The appearance of matrix version of the hyperbolic harmonic
oscillator Lie  group is  significant for several reasons.
 First of all, the hyperbolic harmonic oscillator algebra plays
an important role in attempts of spectral interpretation of the
Riemann zeta-function (see \cite{BK} and \cite{GLO09}). This
relation is to be the subject of our forthcoming publications. Also
note that classical interpretation of the Gamma-function (a basic
building block of the local Archimedean $L$-factors) as a matrix
element of the affine Lie group ${\rm Aff}(\IR)=\IR\rtimes\IR^*$  (see
e.g. \cite{V}) corroborates our interpretation  of the Hecke-Baxter
operator as a matrix element, since ${\rm Aff}(\IR)$ is a subgroup
of the hyperbolic harmonic oscillator group. In particular, spectral
theory of ${\rm Aff}(\IR)$ naturally fits into
framework of the inverted oscillator group. (see e.g \cite{Ho}).
Next, it is worth noting that the hyperbolic harmonic
oscillator group is also instrumental in the description of the
Gelfand-Kapranov-Zelevinsky hypergeometric functions as matrix
elements \cite{GLO23}, \cite{GLO22}. The proposed representation of
the Hecke-Baxter operator by Whittaker function for extensions of
reductive Lie groups might lead to interesting applications. Thus it
would be fruitful to consider matrix generalization of the
Gelfand-Kapranov-Zelevinsky construction \cite{GKZ} by using matrix
hyperbolic harmonic oscillator group we encounter in this note. This
would provide a kind of general construction for the matrix analogs
of various special functions (compare e.g. with \cite{Ka}). Finally
let us note that the Hecke-Baxter operators are also constructed for
other series of classical Lie groups in \cite{GL}. Their
interpretation along the line of this note would be an interesting
and useful task.

The paper organized as follows. Section 2 contains basic facts on
the $GL_{\ell+1}(\IR)$ Hecke-Baxter operator from \cite{GLO08}. In
Section 3 a simple example of the basic construction for
$GL_{1}(\IR)$ is considered. In  Section 4 the general construction
for $GL_{\ell+1}(\IR)$ is presented. Finally in Section 5 we lift
the functions constructed in Section 4 to a matrix element of the
Heisenberg group extension of $Sp_{2\ell+2}\times Sp_{2\ell+2}$.

{\it  Acknowledgements:} The research of S.V.O. is partially
supported by the Beijing Natural Science Foundation grant IS24004.

\section{The Hecke-Baxter  operator for $GL_{\ell+1}(\IR)$}

In this Section we recall the construction \cite{GLO08} of the
Hecke-Baxter operator for the real group $GL_{\ell+1}(\IR)$. For
pedagogical survey of the  spherical Hecke algebra  basics see e.g.
\cite{ILP}.  In the following, we write
$GL_{\ell+1}:=GL_{\ell+1}(\IR)$ for more concise notation. Let
$O_{\ell+1}$ be the orthogonal subgroup in $GL_{\ell+1}$. The
convolution operation defines an action of compactly supported
smooth functions on $GL_{\ell+1}$ on the function space
$L^1(GL_{\ell+1})$. This results in the commutative convolution
algebra structure on the space of compactly supported smooth
functions given by
 \be\label{CONV}
  (f_1*f_2)(g)=\int_G\,d\mu^H(\tilde{g})\,\,f_1(g\tilde{g}^{-1})\,
  f_2(\tilde{g})\,,
 \ee
Equivalently, \eqref{CONV} provides the action of the function $f_1$
on $f_2\in L^1(GL_{\ell+1})$ via the
  associated integral operator with kernel
 \be\label{INTKER}
  K_{f_1}(g,\tilde{g})=f_1(g\tilde{g}^{-1})\,d\mu^H(\tilde{g})\,.
 \ee
Here  $d\mu^H$ is a Haar measure on $GL_{\ell+1}$, which is unique
up to multiplication by an element of $\IR_+$. Thus
 the multiplication structure \eqref{CONV} depends on the
 choice of the Haar measure.  To get rid of this dependence
one might  consider the space of compactly supported  measures on
$GL_{\ell+1}$ with its natural convolution action on the space
$L^1(GL_{\ell+1})$. To recover \eqref{CONV}, consider the measure on
$GL_{\ell+1}$ of the form
 \be\label{MES1}
  d\mu^f(g)=f(g)\, d\mu^H(g)\,,
 \ee
where  $f$ is a positive smooth compactly supported function.
In the following we  use less invariant formulation
\eqref{CONV} with a fixed choice of the Haar measure, but we comment
on the relevance of a more invariant formulation  in terms of
measures at the end of this Section.

Consider the space of functions on $GL_{\ell+1}$ invariant under the
left and right action of  $O_{\ell+1}$. The algebra structure
\eqref{CONV} descends to the subspace of the
$O_{\ell+1}$-biinvariant functions and defines the structure of
spherical Hecke algebra $\mathcal{H}(GL_{\ell+1},O_{\ell+1})$. We
extend the standard definition slightly by allowing
$O_{\ell+1}$-biinvaraint functions to be in the Schwartz space of
$GL_{\ell+1}$.

Let $B_-\subset GL_{\ell+1}(\IR)$ be the Borel subgroup of lower
triangular matrices,  let $N\subset B_-$ be the maximal
unipotent subgroup and let us identify the Cartan torus $A\subset
B_-$ with the subgroup of diagonal matrices, so that $B_-=N\rtimes
A$. Then for $\g\in\IR^{\ell+1}$, let $\chi_{\g}:\,B_-\to\IC^*$ be
the $B_-$-character, for $b=na\in B_-$, given by
 \be\label{CHAR}
  \chi_{\g}(na)=\prod_{j=1}^{\ell+1}|a_j|^{\imath\g_j-\rho_j},\qquad
  \rho_j=\frac{\ell}{2}+1-j\,.
 \ee
Let $(\pi_{\g},\CV_{\g})$ be the spherical principal series
irreducible $GL_{\ell+1}$-representation  induced from \eqref{CHAR}
realized in the space $\CV_{\g}={\rm
Ind}_{B_-}^{GL_{\ell+1}}(\chi_{\g})$ of $B_-$-equivariant functions:
 \be
  f(bg)=\chi_{\g}(b)\,f(g)\,,\quad b\in B_-\,.
 \ee
Then for $g\in G$, the group action in $\CV_{\g}$ is given by
 \be
  \bigl(\pi_{\g}(g)f\bigl)(h)\,=\,f(hg)\,,\quad f\in\CV_{\g}\,.
 \ee
The spherical principal series representation $(\pi_{\g},\CV_{\g})$
is unitarizable, and there exists an invariant skew-Hermitian
pairing $\<\,,\,\>$ on $\CV_{\g}$.

By the multiplicity one theorem \cite{Sha}, there is a unique (up to
multiplication by an element of $\IC^*$)  $O_{\ell+1}$-invariant
vector $\varphi$ in  $\CV_{\g}$. In the following
 we impose the normalization
condition $\<\varphi,\varphi\>=1$. Then there is a linear continuous
map $F:\,\CV_{\g}\to L^2(G)$ sending every non-zero $v\in\CV_{\g}$
to the matrix element of $(\pi_{\g},\CV_{\g})$:
 \be\label{ME}
  F_{v}(g)=\<\varphi,\pi_\gamma(g)\,v\>\,,\quad g\in G\,.
 \ee
Every $O_{\ell+1}$-biinvariant function $\Phi\in
\CH(GL_{\ell+1},O_{\ell+1})$ provides convolution action on
$F_v(g)$, so every matrix element of the form \eqref{ME} is a
$\Phi$-eigenfunction (by uniqueness of $\varphi$) ,
 \be\label{AC1}
  (\Phi*F_{v})(g)=\Lambda_{\Phi}\,F_{v}(g)\,, \quad g\in
  GL_{\ell+1}\,.
 \ee
Therefore, taking the eigenvalues defines a character of the Hecke
algebra
 \be
  \Lambda\,:\quad
  \CH(GL_{\ell+1},O_{\ell+1})\longrightarrow \IC,\,\qquad \Phi\mapsto\Lambda_{\Phi}\,.
 \ee

In the following we consider two specific  matrix elements of the
given principal series representation $(\pi_{\g},\CV_{\g})$. The
first one is the zonal spherical function:
 \be\label{Sph}
  \Phi_{\g}(g)
  =\<\varphi\,, \pi_{\g}(g)\,\varphi\>\,,
   \quad \Phi_{\g}(\Id)=1\,,
 \ee
where $\Id\in GL_{\ell+1}$ is the identity element. To define the
other matrix element, let us pick a principal character $\psi_{N_+}$
of the maximal unipotent subgroup $N_+\subset GL_{\ell+1}$; the
latter can be identified with the subgroup of upper triangular
unipotent matrices. Then the Jacquet-Whittaker function is given by
the matrix element \cite{Ja},
 \be\label{matel}
  \Psi_{\g}(g)
  =\<\varphi\,, \pi_{\g}(g)\,\psi\>\,,
 \ee
where the  Whittaker
vector $\psi\in\CV_{\g}$ is singled out by the
 following $N_+$-equivariance condition:
 \be\label{equivar}
  \Psi_{\g}(gn)=\psi_{N_+}(n)\,
  \Psi_{\g}(g)\,,\quad n\in N_+\,.
 \ee

As a consequence of \eqref{AC1}, for any $\Phi\in\CH(GL_{\ell+1},O_{\ell+1})$, both
matrix elements \eqref{Sph} and \eqref{matel} are
$\Phi$-eigenfunctions with the same eigenvalue:
 \be\label{phiEigen}
  \bigl(\Phi*\Phi_{\g}\bigr)(g)
  =\Lambda_{\Phi}(\g)\,
  \Phi_{\g}(g),\quad
  \bigl(\Phi*\Psi_{\g}\bigr)(g)
  =\Lambda_{\Phi}(\g)\,
  \Psi_{\g}(g)\,.
 \ee
 Following the Langlands philosophy (see e.g. \cite{ILP},\cite{Kna})
to the principal series
 representation $(\pi_{\g}\,,\CV_\g)$ of $GL_{\ell+1}$ one attaches
the corresponding local Archimedean $L$-factor,
 \be\label{Eigenprop}
  L(s)=\prod_{j=1}^{\ell+1}\,
  \pi^{-\frac{s-\i\gamma_j}{2}}\,
  \Gamma\Big(\frac{s-\i\gamma_j}{2}\Big)\,.
 \ee
One should stress that the particular expression \eqref{Eigenprop}
for the local Archimedean $L$-factor implicitly includes
normalization associated with the global structures of number fields
(captured by the functional equation of the corresponding global
zeta-function). Purely local characteristic of representations over
the field of real numbers gives rise to a more general form of
$L$-factor,
 \be\label{Eigenprop1}
  L(s;c)=\prod_{j=1}^{\ell+1}\,
  (2c)^{\frac{s-\i\gamma_j}{2}}\,
  \Gamma\Big(\frac{s-\i\gamma_j}{2}\Big)\,,
 \ee
depending on a parameter $c\in \IR_+$. Thus  local Archimedean
$L$-factor is better understood as an $\IR_+$-torsor. The expression
\eqref{Eigenprop} arises upon  fixing the section $2c=\pi^{-1}$
of the $\IR_+$-torsor.

  One of the main results of \cite{GLO08} may be formulated as follows:\\
   \noindent
{\it  Let $Q_{s,c}(g)$ be a two-parameter family of
$O_{\ell+1}$-biinvariant functions on $GL_{\ell+1}$  given by
 \be\label{UBO}
  Q_{s,c}(g)=\bigl|\det(g)\bigr|^{s+\frac{l}{2}}\,
  e^{-\frac{1}{2c}\Tr(g^{\top}g)}\,,\quad s\in\IC\,,\quad c\in
  \IR_+\,,
 \ee
where $g^{\top}$ is the standard matrix transposition of $g$. Then
with an appropriate normalization of the Haar measure, the function
$Q_{s,c}\in\CH(GL_{\ell+1},O_{\ell+1})$ acts
 on the spherical vector $\varphi\in\CV_{\g}$ via multiplication by
the corresponding $L$-factor \eqref{Eigenprop}.}\\

One should notice that here we imply a bit
   different normalization conventions for the Haar measure, so that
   \eqref{UBO} differs from that in \cite{GLO08} by the factor
   $2^{\ell+1}$, the size of the intersection $O_{\ell+1}\cap
   A=\{\pm1\}^{\ell+1}$ with the Cartan torus $A\subset GL_{\ell+1}$.
   We also adapt the expressions for a more
   general form of $L$-factor \eqref{Eigenprop1}. This more general
   form  actually appears quite naturally in the related Quantum Field
   Theory considerations \cite{GLO11}.

The structure of $Q_{s,c}$ in  \eqref{UBO} is rather simple: it is a
product of the multiplicative $GL_{\ell+1}$-character $|\det(g)|^s$
and the $O_{\ell+1}$-biinvariant Gaussian function. This suggests
considering the $O_{\ell+1}$-biinvariant Gaussian function  as a
part of the measure on the Lie group,
 \be\label{Gauss}
  d\mu^{\text{Gauss}}(g)\,=\,e^{-\frac{1}{2c}\Tr(g^{\top}g)}\,d\mu^H(g)\,,
 \ee
in concordance with \eqref{MES1}.  Appearance of the Gaussian
measure is tied up with the definition of local Archimedean
$L$-factors and ultimately is imposed by functional equations for
the corresponding global zeta-functions \cite{W}.
The underlying meaning of \eqref{Gauss} is that it is the Archimedean counterpart of
 the characteristic function of the subset $G(\IZ_p)\subset G(\IQ_p)$
arising in the non-Archimedean case. Nevertheless its appearance
still remains rather mysterious.

One approach to  understanding the underlying structure of the
Hecke-Baxter operator \eqref{UBO}  would be to look into the theory
of quantum integrable systems. Recall that the
 operator acting on the Whittaker functions \eqref{matel}
by convolution with the kernel function \eqref{UBO} have already
appeared in the theory of quantum integrable theories of Toda type
as an instance of the Baxter operator \cite{GLO08}. Although this
relation with the theory of quantum integrable systems is very
interesting and fruitful,  it does not explain specific choice of
the $O_{\ell+1}$-biinvariant function \eqref{UBO}. Therefore,
further clarifications  of the appearance of the Gaussian function
in \eqref{UBO} are still required.

Below we make a modest step in this direction by proposing a
description of the kernel \eqref{UBO} in terms of representation
theory of a non-reductive extension of the Lie group $GL_{\ell+1}$.

\section{The $GL_1(\IR)$ Hecke-Baxter  operator via  an extension}

Our objective is to describe the kernel \eqref{UBO} of the
Hecke-Baxter operator as a matrix element of some extension of the
Lie group $GL_{\ell+1}$. In this Section we treat the case of the
group $G=GL_1$. Although it looks like a trivial exercise it allows
for rather straightforward generalization to higher ranks.

We start with the Gelfand pair $K\subset G$ with
 \be
  G=GL_1=\IR^*=\IR-\{0\}\,,\quad K=O_1=\{\pm 1\}\subset
  \IR^*\,.
 \ee
The maximal unipotent subgroup of $GL_1$  is trivial, hence the
 Iwasawa decomposition reduces to
 \be
  G=AK=\IR_+\times \{\pm 1\}\,,\quad\IR_+=\{t\in \IR^*|t>0\}\,.
 \ee
Spherical principal series representation
$(\pi_{\g},\,\CV_{\g}),\,\g\in\IR$ is realized in a one-dimensional
$\IC$-vector space with trivial action of $O_1=\{\pm 1\}$ and the
action of the Lie algebra generator $H\in\gl_1=\Lie(GL_1)$ is given
by
 \be\label{GL1rep}
  \pi_{\g}(H)\cdot v=\i\g\,v\,,\quad
  v\in\CV_{\g}\,.
 \ee
We define a skew-Hermitian  $GL_1$-invariant pairing on $\CV_{\g}$
by fixing a vector $v\in \CV_{\g}$,  and requiring
 \be
  \<v,\,\pi_{\g}(H)\,v\>\,=\,-\<\pi_{\g}(H)\,v,\,v\>\,, \qquad \<v,v\>=1  \,.
  \ee
Vector $v$ is both spherical vector  and the Whittaker vector, hence
in the case of $GL_1$ the matrix elements \eqref{Sph} and
\eqref{matel} coincide and are equal to
 \be\label{GL1W}
  \Phi_{\g}(x)=\<v,\,\pi_{\g}(|x|^H)v\>=|x|^{\i\g}\,\qquad x\in \IR^*\,.
  \ee
The function \eqref{GL1W} also gives the character of
$(\pi_{\g},\,\CV_{\g})$.

The $GL_1$  Hecke-Baxter operator \eqref{UBO} is given by
 \be\label{HBGL1}
  Q_{s,c}(x)\,=\,\,|x|^s\,e^{-\frac{1}{2c}|x|^2},\qquad
  x\in \IR^*\,,\quad s\in \IC\,,\quad c\in \IR_+\,.
  \ee
Note that the function \eqref{HBGL1} is obviously invariant under
the subgroup $O_1=\{\pm1\}$ acting via $t\to \pm t$, and thus
indeed  belongs to the spherical Hecke algebra $\CH(GL_1,O_1)$. We
fix the Haar measure as follows
 \be\label{GL1Haar}
  d\mu^H(y)=\frac{dy}{y}\,, \qquad y\in \IR^*\,.
 \ee
The convolution action of $Q_{s,c}$ on the matrix element \eqref{GL1W}
is easy to calculate with the following result, for ${\rm Re}(s)>0$
 \be
  (Q_{s,c}*\Phi_{\g})(x)\,
  =\!\!\int\limits_{GL_1}\!\!d\mu^H(y)\,\,Q_{s,c}(y)\,\Phi_{\g}(xy^{-1})\,
  =\int\limits_{\IR^*}\frac{dy}{y}\,|y|^s\,e^{-\frac{1}{2c} y^2}\,
  |xy^{-1}|^{\i\g}\\
  =\,2|x|^{\i\g}\!\int\limits_0^{\infty}\frac{dy}{y}\,y^{s-\i\g}\,
  e^{-\frac{1}{2c}y^2}\,
  =\,L(s,c)\,\Phi_{\g}(x)\,,
  \ee
where
\be
  L(s,c)\,
  =\,(2c)^{\frac{s-\i\g}{2}}\,
  \Gamma\Big(\frac{s-\i\g}{2}\Big)\,.
  \ee

Thus the matrix element \eqref{GL1W} is indeed a
$Q_{s,c}$-eigenfunction with the  eigenvalue  equal to the
Archimedean $L$-factor \eqref{Eigenprop} with $\ell=0$ attached to
the one-dimensional principal series representation
$(\pi_{\g},\,\CV_{\g})$ of $GL_1$.

Now we construct an extension of the group $GL_1$ allowing to
describe the Hecke-Baxter kernel \eqref{HBGL1} in representation
theory terms. Let us start with the  Heisenberg algebra $\CH^{(1)}$
generated by $X,Y,C$ satisfying the following defining relations
 \be\label{HEIS}
  [X,Y]=C, \qquad [X,C]=0, \qquad [Y,C]=0\,.
 \ee
Consider the $GL_1$-action of on $\CH^{(1)}$ by automorphisms,  so
that the orthogonal subgroup $O_1\subset GL_1$ acts trivially and
the action of the connected component $\IR_+\subset GL_1$ on
$\CH^{(1)}$ is defined by the action of the generator $H$ of
$\gl_1={\rm Lie}(GL_1)$:
 \be\label{HEIS2}
  [H,X]=X, \qquad
  [H,Y]=-Y,\qquad [H,C]=0\,.
 \ee
Introduce the Lie algebra $\mathfrak{h}\gl_1=\gl_1\ltimes\CH^{(1)}$
generated by $H,X,Y,Z$ and with the defining relations \eqref{HEIS}
and \eqref{HEIS2}. We denote by $\CH GL_1$ a Lie group such that
$\mathfrak{h}\gl_1=\Lie(\CH GL_1)$. The Lie algebra
$\mathfrak{h}\gl_1$ may be called the Lie algebra of hyperbolic
harmonic oscillator by  analogy with a more standard case of the Lie
algebra of harmonic oscillator (see e.g. \cite{K}).

Representation theory of the Heisenberg Lie algebra is very simple.
For a fixed nontrivial action of the central element, the Heisenberg
Lie algebra $\CH^{(1)}$ has a unique irreducible unitary
representation. In turn, irreducible representations of
$\mathfrak{h}\gl_1$ (for invertible image of $C$) are directly
related with irreducible representations of $\CH^{(1)}$ and of
$\mathfrak{gl}_1$. To infer this, note that the universal enveloping
algebra $\CU(\mathfrak{h}\gl_1)$ of $\mathfrak{h}\gl_1$ with the
invertible generator $C$ allows the following factorization:
 \be\label{EQUNIV}
  \CU(\mathfrak{h}\gl_1)\otimes_{\IC[C]}
  \IC[C,C^{-1}]\simeq \CU(\gl_1)\otimes \left(\CU(\CH^{(1)})\otimes_{\IC[C]}
  \IC[C,C^{-1}]\right)\,.
 \ee
Indeed, let the algebra in the r.h.s. be  generated by $\tilde{H}\in
\gl_1$ and $X,Y,C\in\CH^{(1)}$ satisfying the relations
 \be
  [X,Y]=C, \qquad [X,C]=0, \qquad [Y,C]=0\,,\\
  \,[\tilde{H},X]=0, \qquad [\tilde{H},Y]=0, \qquad [\tilde{H},C]=0\,.
 \ee
Then one can notice that these relations become equivalent to the
defining relations \eqref{HEIS} and \eqref{HEIS2} for the algebra in
l.h.s. of \eqref{EQUNIV} via taking
 \be
  \tilde{H}=H-\frac{1}{2}C^{-1}(XY+YX).
 \ee
Therefore, for every irreducible $\mathfrak{h}\gl_1$-representation
$(\pi,\CV)$ with invertible action of the central element $C$,
vector space $\CV$ can be identified with a tensor product of an
irreducible $\gl_1$-representation and an irreducible representation
of the Heisenberg algebra $\CH^{(1)}$ (this tensor product evidently
  does not respect the structure of a
  $\mathfrak{h}\gl_1$-representation).

Let $(\pi_{\kappa,\,c}\,,\CV_{\kappa,\,c})$ be irreducible
$\mathfrak{h}\gl_1$-representation realized in the space $\CS$ of
Schwartz functions in variable $x\in \IR$:
 \be\label{REP}
  \pi_{\kappa,\,c}(Y)=-c\pr_x\,,\quad \pi_{\kappa,\,c}(X)=\imath x\,,
  \quad \pi_{\kappa,\,c}(H)=\frac{1}{2}+\imath \kappa+x\pr_x\,,\\
  \pi_{\kappa,\,c}(C)=\imath c\,,
 \ee
where $c\in\IR_+$ is implied. Consider the contragredient
representation $(\pi^{\vee}_{\kappa,c},\CV^{\vee}_{\kappa,c})$
realized in the space of tempered distributions $\CD$ on $\IR$. We
write down a natural Hermitian pairing in the following  explicit
form
 \be\label{pairing}
  \<\chi_1,\chi_2\>=\!
  \int\limits_{\IR}\!dx\,\,\,\overline{\chi_1(x)}\,\,\,\chi_2(x)\,,
  \ee
  implying that $\ov{\chi_1(x)}\,dx$ represents  a tempered
distribution. Then we have
 \be
  \<\chi_1,\pi_{\kappa,c}(U)\chi_2\>
  =-\<\pi^{\vee}_{\kappa,c}(U)\chi_1,\chi_2\>\,,\quad
  U\in\mathfrak{h}\gl_1\,.
 \ee
Assume that the $\mathfrak{h}\gl_1$-representation
$\pi_{\kappa,\,c}$ integrates to the Lie group action of $\CH GL_1$.

Define the  spherical and  Whittaker vectors, $\phi \in
\CV_{\kappa,c}$ and $\psi \in \CV_{\kappa,c}^{\vee}$, by the
following relations:
 \be
  \pi_{\kappa,c}(X-\imath Y)\cdot \phi=0, \qquad
  \pi^{\vee}_{\kappa,c}(X)\cdot\psi=\imath \psi\,.
 \ee
Explicitly, we have
 \be\label{WW}
  \phi(x)=e^{-\frac{1}{2c}\,x^2},\quad\psi(x)=\delta(x-1)\,,
 \ee
where $\delta(x)$ is  the standard delta-function
 considered as a (tempered) distribution.
 Now for the pairing \eqref{pairing}, introduce
the following matrix element
 \be\label{GLHmatel}
  W_{\kappa,c}(y)=\<\psi\,,\pi_{\kappa,c}(|y|^H)\phi\>\,,\qquad
  y\in GL_1=\IR^*\,,
 \ee
which may be called the $\CH GL_1$-Whittaker function (recall that
$O_1\subset GL_1$ acts trivially on vectors \eqref{WW}). Then
substitution of \eqref{WW} into \eqref{GLHmatel}  results in
 \be
  W_{\kappa,c}(y)\,=\,|y|^{\frac{1}{2}+\imath\kappa}\!\int\limits_{\IR}\!dx\,
  \delta(x-1)\,|y|^{\,x\pr_x}e^{-\frac{1}{2c}\,x^2}
  =|y|^{\frac{1}{2}+\imath \kappa}\,e^{-\frac{1}{2c}\,t^2}\,.
 \ee
Therefore, for $s=1/2+\imath \kappa$ we arrive at the expression
\eqref{HBGL1} for the kernel of the $GL_1$ Hecke-Baxter operator:
 \be
  Q_{s,c}(y)\,=\,|y|^s\,e^{-\frac{1}{2c}\,y^2}\,.
  \ee
To realize the general case of a complex $s$ one should consider
non-unitary $\mathfrak{hgl}$-representations
$(\pi_{\g,c},\CV_{\g,c}),\,\g\in\IC$ (we skip these simple
modifications here). Thus we represent the $GL_1$ Hecke-Baxter
kernel by the $\CH GL_1$-Whittaker function.

\section{The $GL_{\ell+1}(\IR)$ Hecke-Baxter  operator via an extension}

In this Section we generalize considerations of the previous
  Section to the case of the general linear
group  $GL_{\ell+1}$.  We start with defining the Lie group
extension $GL_{\ell+1}\subset \CH GL_{\ell+1}$ analogous to the
extension $GL_{1}\subset \CH GL_{1}$ introduced in the previous
Section.

Let $\CH^{(\ell+1)}$ be the standard Heisenberg Lie algebra
associated with the vector space $\Mat_{\ell+1}$ of real
$(\ell+1)\times(\ell+1)$-matrices. Precisely, we consider the
Heisenberg Lie algebra associated with the symplectic vector space,
 \be\label{symplMAT}
  \IR^{2(\ell+1)^2}\,=\,\Mat_{\ell+1}\oplus\Mat_{\ell+1},\quad
  \Big(\begin{smallmatrix}
  Y\\&\\X
  \end{smallmatrix}\Big)\,\in\IR^{2(\ell+1)^2},\\
  \Omega(X,Y)=\Tr\bigl(dX\wedge dY^{\top}\bigr)\,
  =\,\sum_{i,j=1}^{\ell+1}dX_{ij}\wedge dY_{ij}\,.
 \ee
The corresponding  Lie algebra $\CH^{(\ell+1)}$ is generated by
$X_{ij},\,Y_{ij},\,1\leq i,j\leq\ell+1$ and the central element $C$
satisfying the relations
 \be\label{CRL}
  [X_{ij},X_{kl}]=0, \qquad
  [Y_{ij},Y_{kl}]=0\,,\qquad
  [X_{ij},Y_{kl}]=C\delta_{ik}\delta_{jl}\,.
 \ee
Introduce the semi-direct sum $\mathfrak{h}\gl_{\ell+1}$ of the Lie
algebra $\gl_{\ell+1}\oplus\gl_{\ell+1}$ and the Heisenberg Lie
algebra $\CH^{(\ell+1)}$. Let $\{e_{ij}^{L}\,,1\leq
i,j\leq(\ell+1)\}$ be the standard set of generators of the one copy
of $\gl_{\ell+1}$ and $\{e_{ij}^{R}\,,1\leq i,j\leq(\ell+1)\}$ be
the standard set of generators of the other copy of $\gl_{\ell+1}$.
Then $\mathfrak{h}\gl_{\ell+1}$ is defined by the relations
\eqref{CRL} and
 \be\label{OA}
  [e^{L,R}_{ij},e^{L,R}_{kl}]=\delta_{jk}e^{L,R}_{il}
  -\delta_{il}e^{L,R}_{kj}, \quad
  \,[e^L_{ij},e^R_{kl}]=0, \quad
  [e^{L,R}_{ij},C]=0\,,
 \ee
 \be\label{hglL}
  [e^L_{ij},X_{kl}]=-\delta_{ik}X_{jl}, \qquad
  [e^L_{ij},Y_{kl}]=\delta_{jk}Y_{il}\,.
 \ee
 \be\label{hglR}
  [e^R_{ij},X_{kl}]=\delta_{jl}X_{ki}, \qquad
  [e^R_{ij},Y_{kl}]=-\delta_{il}Y_{kj}\,.
 \ee
The Lie algebra $\mathfrak{h}\gl_{\ell+1}$ might be considered as a
special kind of  matrix generalization of the hyperbolic harmonic
oscillator algebra. We denote by $\CH GL_{\ell+1}$ a Lie group such
that $\Lie\bigl(\CH GL_{\ell+1}\bigr)=\mathfrak{h}\gl_{\ell+1}$.

Irreducible representations  of the extension
$\mathfrak{h}\gl_{\ell+1}$ of $\gl_{\ell+1}\oplus\gl_{\ell+1}$ may
be constructed using the standard Mackey technique (see e.g
\cite{K}), but in our case direct methods are more appropriate.

\begin{lem}\label{LEM}
For the universal enveloping Lie algebra of
$\mathfrak{h}\gl_{\ell+1}$ with the invertible generator $C$, the
following factorization holds:
 \be\label{UEFACT}
  \CU(\mathfrak{h}\gl_{\ell+1})\otimes_{\IC[C]}
  \IC[C,C^{-1}]\\
  \simeq\,\CU\bigl(\gl_{\ell+1}
  \oplus\gl_{\ell+1}\bigr)
  \otimes \left(\CU(\CH^{(\ell+1)})
  \otimes_{\IC[C]}  \IC[C,C^{-1}]\right)\,.
 \ee
Furthermore, any  irreducible
$\mathfrak{h}\gl_{\ell+1}$-representation with an invertible action
of the central element $C$ is isomorphic to a representation
$(\pi_{\g_L,\g_R,c},\CV_{\g_L,\g_R,c}),\,\g_L,\g_R\in\IC^{\ell+1},c\in\IR_+$,
such that the underlying vector space $\CV_{\g_L,\g_R,c}$ may be
identified (but not as a representation) with
$\CV_{\g_L,\g_R}\otimes\CW_c$, where
$(\pi_{\g_L,\g_R},\,\CV_{\g_L,\g_R})\simeq(\pi_{\g_L},\,
\CV_{\g_L})\otimes(\pi_{\g_R},\,\CV_{\g_R})$ is an irreducible
$\gl_{\ell+1}\oplus\gl_{\ell+1}$-representation and $\CW_c$ is the
irreducible representation of the Heisenberg Lie algebra
$\CH^{(\ell+1)}$.
\end{lem}

\proof Let us start with explicit construction of the identification
\eqref{UEFACT}. Explicitly  the  algebra in the r.h.s. of \eqref{UEFACT} is
generated by $X_{ij},Y_{ij},C$ and $\hat{e}^{L,R}_{ij}$ subjected to
\eqref{CRL} and the following relations:
 \be
  [\hat{e}^{L,R}_{ij},\hat{e}^{L,R}_{kl}]=
  \delta_{jk}\hat{e}^{L,R}_{il}-\delta_{il}\hat{e}^{L,R}_{kj},\quad
  [\hat{e}^{L}_{ij},\hat{e}^{R}_{kl}]=0,\quad
  [\hat{e}^{L,R}_{ij},C]=0\,,\\
  \,[\hat{e}^{L,R}_{ij},X_{kl}]=0, \quad
  [\hat{e}^{L,R}_{ij},Y_{kl}]=0\,,\quad1\leq i,j,k,l\leq(\ell+1)\,,
 \ee
where the generators $\hat{e}^{L,R}_{ij}$ are given by
 \be\label{GLosc}
  \hat{e}^L_{ij}\,=\,e^L_{ij}\,
  -\,C^{-1}\sum_{k=1}^{\ell+1}X_{jk}Y_{ik}\,,
  \quad
  \hat{e}^R_{ij}\,=\,e^R_{ij}\,
  +\,C^{-1}\sum_{k=1}^{\ell+1}X_{ki}Y_{kj}\,\,.
  \ee
The above expressions provide the isomorphism \eqref{UEFACT}. Taking
into account \eqref{UEFACT}, for every irreducible
$\mathfrak{h}\gl_{\ell+1}$-representation $(\pi,\CV)$ with
invertible image of $C$, vector space $\CV$ may be identified with a
tensor product of a $\gl_{\ell+1}\oplus\gl_{\ell+1}$-irreducible
representation and an irreducible representations of the Heisenberg
algebra $\CH^{(\ell+1)}$. By the Stone-von Neumann theorem, for a
fixed invertible image of the central generator $C$, the Lie algebra
$\CH^{(\ell+1)}$ has a unique irreducible representation, therefore
irreducible representations of the Lie algebra
$\mathfrak{h}\gl_{\ell+1}$ (again with invertible image of $C$) are
in one to one correspondence with irreducible representations of
$\gl_{\ell+1}\oplus\gl_{\ell+1}$. $\Box$

Consider the $\CH GL_{\ell+1}$-representation
$(\pi_{\g_L,\g_R,c},\CV_{\g_L,\g_R,c})$ obtained by lifting (in the
sense of Lemma \ref{LEM}) the spherical principal series
$GL_{\ell+1}\times GL_{\ell+1}$-representation
$\CV_{\g_L,\g_R}=\CV_{\g_L}\otimes\CV_{\g_R}$ via representation
$\CW_c,\,c\in\IR_+$ of the Heisenberg group with the central element
$C$ taking value $\imath c\in\imath\IR_+$. By Lemma \ref{LEM}, the
underlying
  vector space $\CV_{\g_L,\g_R,c}$ may be identified with the
tensor product $\CV_{\g_L}\otimes\CV_{\g_R}\otimes \CW_c$. We
consider a realization of the Heisenberg group representation
$\CW_c$ in the space $\CS(\Mat_{\ell+1})$ of Schwartz functions on
the space $\Mat_{\ell+1}$ of real $(\ell+1)\times (\ell+1)$-matrices
$x=\|x_{ij}\|$. The $\mathfrak{h}\gl_{\ell+1}$-action on smooth
vectors in $\CV_{\g_L,\g_R,c}$ is defined by, for $1\leq
i,j\leq(\ell+1)$,
 \be\label{HGLrep}
  \pi_{\g_L,\g_R,c}(X_{ij})=\imath x_{ij}\,, \quad
  \pi_{\g_L,\g_R,c}(Y_{ij})=-c\frac{\pr}{\pr x_{ij}}\,,\quad
  \pi_{\g_L,\g_R,c}(C)=\imath c\,,\\
  \pi_{\g_L,\g_R,c}(e^L_{ij})
  =-\sum_{k=1}^{\ell+1}
  x_{jk}\frac{\pr}{\pr x_{ik}}-\frac{\ell+1}{2}\delta_{ij}
  +\pi_{\g_L}(\hat{e}^L_{ij})\,,\\
  \pi_{\g_L,\g_R,c}(e^R_{ij})
  =\sum_{k=1}^{\ell+1}x_{ki}\frac{\pr}{\pr x_{kj}}
  +\frac{\ell+1}{2}\delta_{ij}
  +\pi_{\g_R}(\hat{e}^R_{ij})\,,
 \ee
where $\pi_{\g_L,\g_R}(\hat{e}^{L,R}_{ij})$ denotes the action of
the $\gl_{\ell+1}\oplus\gl_{\ell+1}$-generators in the spherical
principal series representation
$\CV_{\g_L,\g_R}=\CV_{\g_L}\otimes\CV_{\g_R}$. Let
$(\pi^{\vee}_{\g_L,\g_R,c},\CV_{\g_L,\g_R,c}^{\vee})$ be the
contragredient representation in the dual space
$\CV_{\g_L,\g_R,c}^{\vee}\simeq\CV^{\vee}_{\g_L}\otimes\CV^{\vee}_{\g_R}\otimes
\CW_c^{\vee}$, where the space $\CW_c^{\vee}$ is given by the space
$\CD(\Mat_{\ell+1})$ of tempered distributions on $\Mat_{\ell+1}$.
We write down a natural pairing in the following form
 \be\label{pairing1}
  \<\chi_1,\chi_2\>\,=\!\!\!
  \int\limits_{\Mat_{\ell+1}}\!\!\!\!\!d\mu(x)\,\,\,
  \overline{\chi_1(x)}\,\,\chi_2(x)\,,\quad
  d\mu(x)=\prod_{i,j=1}^{\ell+1}dx_{ij}\,,
 \ee
so that $\ov{\chi_1(x)}\,d\mu(x)$ represents  a tempered
distribution on $\Mat_{\ell+1}$. According to the definition of
contragredient representation,
 \be
  \<\chi_1,\pi_{\g_L,\g_R,c}(U)\chi_2\>
  =-\<\pi^{\vee}_{\g_L,\g_R,c}(U)\chi_1,\chi_2\>\,,\quad
  U\in\mathfrak{h}\gl_{\ell+1}\,.
 \ee
Assume the $\gl_{\ell+1}\oplus \gl_{\ell+1}$-action via
$\pi_{\g_L,\g_R,c}(e_{ij})$ integrates to the $GL_{\ell+1}\times
GL_{\ell+1}$-action which supplies $\CV_{\g_L,\g_R,c}$ with a
structure of $\CH GL_{\ell+1}$-module. Recall the decomposition
 \be\label{HGL}
  \mathfrak{h}\gl_{\ell+1}\,
  =\,\bigl(\gl_{\ell+1}\oplus
  \gl_{\ell+1}\bigr)\ltimes\CH^{(\ell+1)}\,.
 \ee
For the subgroup $O_{\ell+1}\subset GL_{\ell+1}$ of orthogonal
matrices,
  let $K^{L,R}_{ij}=e^{L,R}_{ij}-e^{L,R}_{ji},\,i<j$ be
generators of the two copies of spherical subalgebras
$\mathfrak{so}_{\ell+1}=\Lie(O_{\ell+1})\subset\mathfrak{gl}_{\ell+1}$.
Then by \eqref{hglL} and \eqref{hglR} the following relations hold:
 \be
  [e^L_{ij}-e^L_{ji},X_{kl}]=\delta_{jk}X_{il}-\delta_{ik}X_{jl},\quad
  [e^L_{ij}-e^L_{ji}, Y_{kl}]=-\delta_{ik}Y_{jl}+\delta_{jk}Y_{il},\\
  \,[e^R_{ij}-e^R_{ji},X_{kl}]=-\delta_{il}X_{kj}+\delta_{jl}X_{ki},\quad
  [e^R_{ij}-e^R_{ji},Y_{kl}]=\delta_{jl}Y_{ki}-\delta_{il}Y_{kj}\,.
 \ee
As a consequence we have for arbitrary $\alpha\in\IC$,
 \be\label{HGLsubalg}
  \bigl[e^L_{ij}-e^L_{ji},(X_{kl}+\alpha Y_{kl})\bigr]
  =\delta_{jk}(X_{il}+\alpha Y_{il})
  -\delta_{ik}(X_{jl}+\alpha Y_{jl})\,,\\
  \bigl[e^R_{ij}-e^R_{ji}, (X_{kl}+\alpha Y_{kl})\bigr]
  =-\delta_{il}(X_{kj}+\alpha Y_{kj})
  +\delta_{jl}(X_{ki}+\alpha Y_{ki})\,.
 \ee

Next, given the $\CH GL_{\ell+1}$-representation
$(\pi_{\g_L,\g_R,c},\CV_{\g_L,\g_R,c})$, define the spherical vector
$\phi \in \CV_{\g_L,\g_R,c}$ by the following relations,
 \be\label{WVR}
  \pi_{\g_L,\g_R,c}(X_{ij})\cdot\phi\,
  =\,\imath\,\pi_{\g_L,\g_R,c}(Y_{ij})\cdot\phi, \qquad
  \pi_{\g_L,\g_R,c}(K^{L,R}_{ij})\cdot \phi\,=\,0\,,
 \ee
that are consistent due to \eqref{HGLsubalg}. The equations
\eqref{WVR} allow the following solution:
 \be\label{LRsph}
  \phi(x)\,=\,\varphi_{\g_L}\otimes\varphi_{\g_R}\otimes\,
  e^{-\frac{1}{2c}\,\Tr(x^\top x)}\,,\qquad
  \varphi_{\g_{L,R}}\in\CV_{\g_{L,R}}\,,
 \ee
where $\varphi_{\g_{L,R}}$ are the $O_{\ell+1}$-spherical vectors in
the $GL_{\ell+1}$-representations $(\pi_{\g_{L,R}},\CV_{\g_{L,R}})$
normalized by $\<\varphi_{\g_{L,R}},\,\varphi_{\g_{L,R}}\>=1$.

For the contragredient representation
$(\pi^{\vee}_{\g_L,\g_R,c},\CV_{\g_L,\g_R,c}^{\vee})$, introduce the
Whittaker vector $\psi \in \CV^{\vee}_{\g_L,\g_R,c}$ defined by
 \be\label{WVL}
  \pi^{\vee}_{\g_L,\g_R,c}(X_{ij})\cdot \psi\,
  =\,\imath \delta_{ij}\psi, \qquad
  \pi^{\vee}_{\g_L,\g_R,c}(e^L_{ij}+e^R_{ij})\cdot \psi\,=\,0\,,
  \ee
so that the consistency follows by \eqref{hglL},\eqref{hglR}.
Finding a non-trivial solution to the equations \eqref{WVL} on the
Whittaker vector $\psi\in \CV^{\vee}_{\g_L,\g_R,c}$ requires
imposing $\g_L+\g_R=0$ so that $\CV_{\g_L}\simeq\CV_{\g_R}^{\vee}$.
In the following we denote $\g_R=\g$ and keep the notation
$(\pi_{\g,c},\CV_{\g,c})$ for the representation
$(\pi_{-\g,\g,c},\CV_{-\g,\g,c})$. Then under these assumptions, the
equations \eqref{WVL} allow the following solution,
 \be\label{Whitvec}
  \psi(x)\,=\,\Id_{\End(\CV_{\g})}
  \otimes\prod_{i,j=1}^{\ell+1}\delta(x_{ij}-\delta_{ij})\,,\qquad
  \Id_{\End(\CV_{\g})}\in\End(\CV_{\g})
  \,,
 \ee
where $\Id_{\End(\CV_{\g})}$ is the canonical unit element of the
algebra $\End(\CV_{\g})\simeq\CV_{\g}^{\vee}\otimes\CV_{\g}$ given
by the identity endomorphism. Note that the former relations in
\eqref{WVL} and the delta-factor
$\prod\limits_{i,j=1}^{\ell+1}\delta(x_{ij}-\delta_{ij})$ in
 \eqref{Whitvec} are invariant under the adjoint action of the
diagonal subgroup $GL_{\ell+1}\subset GL_{\ell+1}\times
GL_{\ell+1}$.

Now we introduce the $\CH GL_{\ell+1}$-matrix element given by the
pairing \eqref{pairing1} with  vectors $\phi \in \CV_{\g,c}$ and
$\psi\in \CV^{\vee}_{\g,c}$:
 \be\label{GUBO0}
  W_{\gamma,c}(g_1,g_2)\,
  =\,\<\psi,\,\pi_{\gamma,c}(g_1\otimes g_2)\,\phi\>\,,\quad \
  g_1\otimes g_2\in GL_{\ell+1}\times GL_{\ell+1}\,.
 \ee
Then a direct calculation gives
 \be\label{GUBO}
  W_{\gamma,c}(g_1,g_2)\,
  =\!\int\limits_{\Mat_{\ell+1}}\!\!\!dx
  \prod_{i,j=1}^{\ell+1}\delta(x_{ij}-\delta_{ij})\,
  \Phi_{\g}(g_1^{-1}g_2)\,\,e^{-\frac{1}{2c}\Tr[(g_1^{-1}xg_2)^{\top}g_1^{-1}xg_2]}\\
  =\,\Phi^{\gl_{\ell+1}}_{\gamma}(g_1^{-1}g_2)\,\,
  e^{-\frac{1}{2c}\Tr[(g_1g_2^{-1})^{\top}g_1g_2^{-1}]}\,,
 \ee
where $\Phi^{\gl_{\ell+1}}_{\gamma}(g)$ is the zonal spherical
function of the principal series $\gl_{\ell+1}$-representation
$\CV_{\gamma}$ defined in \eqref{Sph}. Actually the function
\eqref{GUBO} on $GL_{\ell+1}\times GL_{\ell+1}$ is an integral
kernel of the convolution with the function
 \be\label{GUBO1}
  W_{\gamma,c}(g)
  =\<\psi,\,\pi_{\gamma_1,\g_2,c}(\Id\otimes g)\,\phi\>=
  \Phi^{\gl_{\ell+1}}_{\gamma}(g)\,\,
  e^{-\frac{1}{2c}\Tr g^{\top}g}\,,
 \ee
with respect to the Haar measure on $GL_{\ell+1}$.

\begin{te}
For $\kappa\in\IR$, let $(\pi_{\kappa},\CV_\kappa)$  be the
one-dimensional representation of $GL_{\ell+1}\times GL_{\ell+1}$
given by
 \be\label{1dimRP}
  \pi_\kappa\,:\quad
  (g_1,g_2)\longmapsto\bigl|\det(g_1^{-1}g_2)\bigr|^{\imath \kappa}\,,
 \ee
and let $\CW_c,\,c\in \IR_+$ be the irreducible
$\CH^{(\ell+1)}$-representation with  the image of central element
being $\imath c$. Let $(\pi_{\kappa,c}\,,\CV_{\kappa,c})$ be the
irreducible $\CH GL_{\ell+1}$-module with the representation space
$\CV_{\kappa,c}\simeq\CV_\kappa^{\vee}\otimes
\CV_\kappa\otimes\CW_c$ according to Lemma \ref{LEM}. Then the
$GL_{\ell+1}$ Hecke-Baxter operator defined in \eqref{UBO},
 \be\label{UBO12}
  Q_{s,c}(g)=\bigl|\det(g)\bigr|^{s+\frac{s}{2}}\,\,
  e^{-\frac{1}{2c}\Tr(g^{\top}g)}\,,
 \ee
allows the following $\CH GL_{\ell+1}$  matrix element
representation:
 \be
  Q_{s,c}(g)=\<\psi,\,\pi_{\kappa,c}(1,g)\,\phi\>\,,\qquad
  s=\imath\kappa+\frac{1}{2}\,,
  \quad g\in
  GL_{\ell+1}\,.
 \ee
where $\phi$ is given by \eqref{LRsph} and $\psi$ is given by
\eqref{Whitvec}.
\end{te}

\proof Given the $\CH GL_{\ell+1}$-representation
$(\pi_{\g_L,\g_R,c},\CV_{\g_L,\g_R,c})$ labeled by
$\g_L,\g_R\in\IR^{\ell+1},c\in\IR_+$, by Lemma \ref{LEM}, the
representation space $\CV_{\g_L,\g_R,c}$ may be realized (as a
vector space) by the tensor product
$\CV_{\g_L,\g_R,c}=\CV_{\g_L,\g_R}\otimes\CW_c$ of irreducible
$GL_{\ell+1}\times GL_{\ell+1}$-representation
$(\pi_{\g_L,\g_R},\CV_{\g_L,\g_R})$ with
$\CV_{\g_L,\g_R}=\CV_{\g_L}\otimes\CV_{\g_R}$ and irreducible
$\CH^{(\ell+1)}$-representation $\CW_c$. Specifying the
$GL_{\ell+1}\times GL_{\ell+1}$-representation
$\CV_{\g_L,\g_R}=\CV_{\g_L}\otimes\CV_{\g_R}$ to be a
one-dimensional representation \eqref{1dimRP} the corresponding
spherical function \eqref{Sph} equals
 \be
  \Phi^{\gl_{\ell+1}}_{\kappa}(g)=|\det(g)|^{\imath \kappa}\,.
 \ee
Substituting this into \eqref{GUBO} results in
 \be
  W_{\kappa,c}(g)\,
  =\,|\det g|^{\imath \kappa+\frac{\ell+1}{2}}\,
  e^{-\frac{1}{2c} \Tr(g^{\top}g)}\,.
 \ee
Thus we have reconstructed the Hecke-Baxter operator \eqref{UBO} for
$s=\imath \kappa+1/2$  as a special case of $\CH
GL_{\ell+1}$-Whittaker function. $\Box$

Actually, the matrix elements \eqref{GUBO0}, \eqref{GUBO1} provide a
natural generalization of the Hecke-Baxter operator \eqref{UBO12}
with the integral kernel given by the Whittaker function
\eqref{GUBO}. It would be interesting to find the corresponding
eigenvalues for convolution action of the general integral operators
\eqref{GUBO0} on $GL_{\ell+1}$-spherical functions.  While leaving
this for another occasion, we only notice here a trivial fact that
the action of the generalized Hecke-Baxter operator \eqref{GUBO} in
a one-dimensional $GL_{\ell+1}$-representation reduces to
multiplication by the local Archimedean $L$-factor
\eqref{Eigenprop}.

\section{Lifting to symplectic Lie group matrix elements}

In this Section we construct a lifting of the $\CH
GL_{\ell+1}$-Whittaker function considered in the previous Section
to a matrix element of the larger Lie group $\CH Sp_{2\ell+2}$ given
by the semidirect product of $Sp_{2\ell+2}\times Sp_{2\ell+2}$
and the Heisenberg group accosiated with the matrix space
$\Mat_{\ell+1}$. Let us start with description of the corresponding
Lie algebra $\CH\mathfrak{sp}_{2\ell+2}={\rm Lie}(\CH
Sp_{2\ell+2})$.

Recall that the symplectic Lie algebra $\mathfrak{sp}_{2\ell+2}$
allows the following embedding:
 \be\label{SPLiealg}
  \ssp_{2\ell+2}\,=\,\Big\{
  \Big(\begin{smallmatrix}
  A&&B\\&&\\C&&-A^{\top}
  \end{smallmatrix}\Big)\,:\quad B^{\top}=B,\quad C^{\top}=C\Big\}\,
  \subset\,\gl_{2\ell+2}\,,
 \ee
with the block-diagonal subalgebra
$\gl_{\ell+1}=\diag\{A,\,-A^{\top}\}\subset\ssp_{2\ell+2}$. Denote
by $a_{ij},\,1\leq i,j\leq(\ell+1)$ and $b_{ij}, c_{ij},\,1\leq
i\leq j\leq(\ell+1)$ the corresponding generators of
$\mathfrak{sp}_{2\ell+2}$. In the following we consider two copies
of $\mathfrak{sp}_{2\ell+2}$ with the two sets of generators:
$a^L_{ij}, b^L_{ij}, c^L_{ij}$ and $a^R_{ij}, b^R_{ij}, c^R_{ij}$.
The generators $a^{L,R}_{ij}\,,b^L_{ij},c^L_{ij}\,,1\leq i\leq j\leq
\ell+1$ of the Lie algebra $\ssp_{2\ell+2}\oplus\ssp_{2\ell+2}$
satisfy the following commutation relations:
 \be\label{SPrel}
  [b^{L,R}_{ij},b^{L,R}_{kl}]=[c^{L,R}_{ij},c^{L,R}_{kl}]=0, \quad
  [a^{L,R}_{ij},a^{L,R}_{kl}]= \delta_{jk}a^{L,R}_{il}-\delta_{il}a^{L,R}_{kj}, \\
  \,[b^{L,R}_{ij},c^{L,R}_{kl}]=\delta_{ik}a^{L,R}_{jl}+\delta_{jl}a^{L,R}_{ik}
  +\delta_{jk}a^{L,R}_{il}+\delta_{il}a^{L,R}_{jk},\\
  \,[a^{L,R}_{ij},b^{L,R}_{kl}]=\delta_{jk}b^{L,R}_{il}+\delta_{jl}b^{L,R}_{ik}, \quad
  [a^{L,R}_{ij},c^{L,R}_{kl}]=-\delta_{ik}c^{L,R}_{jl}-\delta_{il}c^{L,R}_{jk}\,.
 \ee

Now we introduce the semidirect sum
$\fh\ssp_{2\ell+2}=\bigl(\ssp_{2\ell+2}\oplus\ssp_{2\ell+2}\bigr)
\ltimes\CH^{(\ell+1)}$ as follows. The Lie algebra
$\fh\ssp_{2\ell+2}$ is generated by
$a_{ij}^{L,R},X_{ij},Y_{ij}\,,1\leq i,j\leq\ell+1$, and
$b_{ij}^{L,R},c_{ij}^{L,R},\,1\leq j$ and the central element $C$
satisfying the commutation relations \eqref{SPrel} and
\eqref{CRL},\eqref{hglL},\eqref{hglR} with the following additional
relations:
 \be\label{HSPrelL}
  [b^L_{ij},Y_{kl}]=0,\quad
  [b^L_{ij},X_{kl}]=\delta_{jk}Y_{il}+\delta_{ik}Y_{jl},\\
  \,[c^L_{ij},X_{kl}]=0,\quad
  [c^L_{ij},Y_{kl}]=\delta_{jk}X_{il}+\delta_{ik}X_{jl}\,,
 \ee
and
 \be\label{HSPrelR}
  [b^R_{ij},X_{kl}]=0,\quad
  [b^R_{ij},Y_{kl}]=-\delta_{il}X_{kj}-\delta_{jl}X_{ki},\\
  \,[c^R_{ij},Y_{kl}]=0,\quad
  [c^R_{ij},X_{kl}]=-\delta_{il}Y_{kj}-\delta_{jl}Y_{ki}\,.
 \ee
It is easy to check directly that thus defined action of the Lie
algebra  of $ \ssp_{2\ell+2}\oplus\ssp_{2\ell+2}$ is compatible with
the structure of the Heisenberg Lie algebra
 \be\label{CRL1}
  [X_{ij},X_{kl}]=0, \qquad
  [Y_{ij},Y_{kl}]=0\,,\qquad
  [X_{ij},Y_{kl}]=C\delta_{ik}\delta_{jl}\,.
 \ee
Using the injective homomorphism,
 \be\label{GL2SP}
  \gl_{\ell+1}\oplus\gl_{\ell+1}\,
  \longrightarrow\,
  \ssp_{2\ell+2}\oplus\ssp_{2\ell+2}\,,\\
  (e^L_{ij}\,,e^R_{ij})\longmapsto(a^L_{ij}\,,a^R_{ij})\,,\quad1\leq
  i,j\leq\ell+1\,,
  \ee
we obtain the embedding of the Lie algebra \eqref{HGL} into
$\hsp_{2\ell+2}$:
 \be\label{hglSUB}
  \mathfrak{h}\gl_{\ell+1}\,
  =\,\bigl(\gl_{\ell+1}\oplus
  \gl_{\ell+1}\bigr)\ltimes\CH^{(\ell+1)}\,
  \subset\,\fh\ssp_{2\ell+2}\,.
 \ee
In the following, we consider representations of the Lie algebra
$\hsp_{2\ell+2}$. Similarly to the $\hgl_{\ell+1}$-representations
$(\pi_{\g_L,\g_R,c},\CV_{\g_L,\g_R,c})$ considered in the previous
Section, the $\hsp_{2\ell+2}$-representations are labeled by the
weights $\g=(\g_L,\g_R),\,\g_{L,R}\in\IC^{\ell+1}$ of the Lie
algebra $\ssp_{2\ell+2}\oplus\ssp_{2\ell+2}$ and by the value
$\imath c$ of $C\in\CH^{(\ell+1)}$. Note that  we keep the same
notations for the labels of representations as in the case of
$\gl_{\ell+}\oplus\gl_{\ell+1}$ by taking into account the embedding
\eqref{hglSUB}.
\begin{lem}\label{LEM1}
For the universal enveloping
Lie algebra of $\mathfrak{hsp}_{2\ell+2}$ with the invertible
generator $C$, the following factorization holds
 \be\label{USPfact}
  \CU(\mathfrak{h}\ssp_{2\ell+2})\otimes_{\IC[C]}
  \IC[C,C^{-1}]\\
  \simeq\, \CU\bigl(\ssp_{2\ell+2}
  \oplus\ssp_{2\ell+2}\bigr)
  \otimes \left(\CU(\CH^{(\ell+1)})
  \otimes_{\IC[C]}  \IC[C,C^{-1}]\right)\,.
 \ee
Furthermore, any  irreducible $\mathfrak{h
sp}_{2\ell+2}$-representation with an invertible action of the
central element $C$ is isomorphic to a representation
$(\pi_{\g_L,\g_R,c},\CV_{\g_L,\g_R,c})$ labeled by
$\g_L,\g_R\in\IC^{\ell+1},c\in\IR_+$, such that underlying vector
space $\CV_{\g_L,\g_R,c}$ may be identified (but not as a
representation) with $\CV_{\g_L,\g_R}\otimes\CW_c$, where
$\CV_{\g_L,\g_R}\simeq \CV_{\g_L}\otimes \CV_{\g_R}$ is an
irreducible $\ssp_{2\ell+2}\oplus\ssp_{2\ell+2}$-representation and
$\CW_c$ is the irreducible representation of the Heisenberg algebra
$\CH^{(\ell+1)}$.
\end{lem}

\proof The proof of \eqref{USPfact} is the straightforward
generalization of the that for Lemma \ref{LEM}. Let the algebra
$\CU\bigl(\ssp_{2\ell+2}\oplus\ssp_{2\ell+2}\bigr)$ on the r.h.s. of
\eqref{USPfact}  is generated by the generators
$\hat{a}_{ij}^{L,R},\,1\leq i,j\leq\ell+1$ and
$\hat{b}_{ij}^{L,R},\hat{c}_{ij}^{L,R},\,1\leq j$. The generators
$\hat{a}^{L,R}_{ij}$ are identified with $\hat{e}^{L,R}_{ij}$
\eqref{GLosc} due to \eqref{GL2SP}, and the other generators, for
$i\leq j$, are given by
 \be
  \hat{b}_{ij}^L
  =b_{ij}^L\,+\,C^{-1}\sum_{k=1}^{\ell+1}Y_{jk}Y_{ik}\,,\quad
  \hat{b}_{ij}^R
  =b_{ij}^R\,-\,C^{-1}\sum_{k=1}^{\ell+1}X_{ki}X_{kj}\,,\\
  \hat{c}_{ij}^L
  =c_{ij}^L\,-\,C^{-1}\sum_{k=1}^{\ell+1}X_{jk}X_{ik}\,,\quad
  \hat{c}_{ij}^R
  =c_{ij}^R\,+\,C^{-1}\sum_{k=1}^{\ell+1}Y_{ki}Y_{kj}\,.
 \ee
As a consequence of \eqref{USPfact}  every irreducible
$\mathfrak{h}\ssp_{2\ell+2}$-representation with an invertible
action of the central element $C$ is equivalent to a representation
$(\pi_{\g_L,\g_R,c},\CV_{\g_L,\g_R,c}),\,\g_L,\g_R\in\IC^{\ell+1},c\in\IR_+$,
such that the vector space $\CV_{\g_L,\g_R,c}$ may be identified
(not as a representation) with $\CV_{\g_L,\g_R}\otimes\CW_c$, where
$\CV_{\g_L,\g_R}\simeq \CV_{\g_L}\otimes \CV_{\g_R}$ is an
irreducible $\ssp_{2\ell+2}\oplus\ssp_{2\ell+2}$-representation and
$\CW_c$ is the irreducible representation of the Heisenberg algebra
$\CH^{(\ell+1)}$. $\Box$

Consider the $\CH Sp_{2\ell+2}$-representation
$(\pi_{\g_L,\g_R,c},\CV_{\g_L,\g_R,c})$ obtained by lifting the
spherical principal series $Sp_{2\ell+2}\times
Sp_{2\ell+2}$-representation in the space
$\CV_{\g_L,\g_R}=\CV_{\g_L}\otimes\CV_{\g_R}$ via the representation
$\CW_c,\,c\in\IR_+$ of the Heisenberg group with the central element
$C$ taking value $\imath c\in\imath\IR_+$. We consider a realization
of the Heisenberg group representation $\CW_c$ in the space
$\CS(\Mat_{\ell+1})$ of Schwartz functions on the matrix space
$\Mat_{\ell+1}$, so that the subalgebra
$\mathfrak{h}\gl_{\ell+1}\subset\fh\ssp_{2\ell+2}$ acts on smooth
vectors in $\CV_{\g_L,\g_R,c}$ via \eqref{HGLrep} and the action of
other generators $b^{L,R}_{ij},c^{L,R}_{ij},\,i\leq j$ is given by
 \be
  b_{ij}^L
  =\imath c\sum_{k=1}^{\ell+1}\frac{\pr}{\pr x_{ik}\pr x_{jk}}\,
  +\,\pi_{\g_L,\g_R}(\hat{b}^L_{ij})\,,\quad
  b_{ij}^R
  =\frac{\imath}{c}\sum_{k=1}^{\ell+1}x_{ki}x_{kj}\,
  +\,\pi_{\g_L,\g_R}(\hat{b}^R_{ij})\,,\\
  c_{ij}^L
  =\frac{\imath}{c}\sum_{k=1}^{\ell+1}x_{ik}x_{jk}\,
  +\,\pi_{\g_L,\g_R}(\hat{c}^L_{ij})\,,\quad
  c_{ij}^R
  =\imath c\sum_{k=1}^{\ell+1}\frac{\pr}{\pr x_{ki}\pr x_{ki}}\,
  +\,\pi_{\g_L,\g_R}(\hat{c}^R_{ij})\,.
 \ee
where
$\pi_{\g_L,\g_R}(\hat{b}^{L,R}_{ij}),\,\pi_{\g_L,\g_R}(\hat{c}^{L,R}_{ij})$
denote the action of $\ssp_{2\ell+2}\oplus\ssp_{2\ell+2}$-generators
in the spherical principal series representation space
$\CV_{\g_L,\g_R}=\CV_{\g_L}\otimes\CV_{\g_R}$. Let
$(\pi^{\vee}_{\g_L,\g_R,c},\CV_{\g_L,\g_R,c}^{\vee})$ be the
contragredient representation in the dual space
$\CV_{\g_L,\g_R,c}^{\vee}\simeq\CV^{\vee}_{\g_L}\otimes\CV^{\vee}_{\g_R}\otimes
\CW_c^{\vee}$, where the space $\CW_c^{\vee}$ is given by the space
$\CD(\Mat_{\ell+1})$ of tempered distributions on $\Mat_{\ell+1}$.
We write down the corresponding  pairing in the following  form
 \be\label{pairing2}
  \<\chi_1,\chi_2\>\,=\!\!\!
  \int\limits_{\Mat_{\ell+1}}\!\!\!\!d\mu(x)\,\,\,
  \overline{\chi_1(x)}\,\,\chi_2(x)\,,\quad
  d\mu(x)=\prod_{i,j=1}^{\ell+1}dx_{ij}\,,
 \ee
so that $\ov{\chi_1(x)}\,d\mu(x)$ represents  a tempered
distribution on $\Mat_{\ell+1}$. The contragredient representation
is defined by
 \be
  \<\chi_1,\pi_{\g_L,\g_R,c}(U)\chi_2\>
  =-\<\pi^{\vee}_{\g_L,\g_R,c}(U)\chi_1,\chi_2\>\,,\quad
  U\in\mathfrak{h}\ssp_{2\ell+2}\,.
 \ee
Assume the $\ssp_{2\ell+2}\oplus\ssp_{2\ell+2}$-action via
$\pi_{\g_L,\g_R,c}$ integrates to the $Sp_{2\ell+2}\times
Sp_{2\ell+2}$-action which supplies $\CV_{\g_L,\g_R,c}$ with a
structure of $\CH Sp_{2\ell+2}$-module.

Next, we construct lifting of the defining relations
\eqref{WVL},\eqref{WVR} for the vectors $\phi,\psi$.

\begin{lem}\label{SPLEM}
  The subspace $\fk\subset\fh\ssp_{2\ell+2}$ spanned by
 \be\label{HSPsphalg}
  K^{L,R}_{ij}\,=\,a^{L,R}_{ij}-a^{L,R}_{ji}\,,\quad1\leq
  i<j\leq\ell+1,\\
  b^{L,R}_{ij}-c^{L,R}_{ij}\,,\quad1\leq i\leq j\leq\ell+1,\\
  X_{ij}-\imath Y_{ij}\,,\quad1\leq i,j\leq\ell+1\,.
 \ee
is a subalgebra isomorphic to the semidirect sum:
 \be
  \fk\,\simeq\,\bigl(\fu_{\ell+1}\oplus\fu_{\ell+1}\bigr)
  \ltimes\CA,\quad
  \CA=\Span\bigl\{X_{ij}-\imath Y_{ij}\,,\,\,
  1\leq i,j\leq\ell+1\bigr\}\,.
 \ee
\end{lem}
\proof The Lie subalgebra of $\fk$ spanned by the elements
 \be
  K^{L,R}_{ij}\,=\,a^{L,R}_{ij}-a^{L,R}_{ji}\,,\quad1\leq
  i<j\leq\ell+1,\\
  b^{L,R}_{ij}-c^{L,R}_{ij}\,,\quad1\leq i\leq j\leq\ell+1
 \ee
is isomorphic to the spherical (maximal compact) subalgebra of
$\ssp_{2\ell+2}\oplus\ssp_{2\ell+2}$ i.e. a direct sum of unitary
Lie algebras:
 \be
  \fu_{\ell+1}\oplus\fu_{\ell+1}
  \subset\ssp_{2\ell+2}\oplus\ssp_{2\ell+2}\,.
 \ee
Consider the subalgebra $\CA=\Span\{X_{ij}-\imath Y_{ij}\,,\,
  1\leq i,j\leq\ell+1\}$, which is commutative:
 \be
  [X_{ij}-\imath Y_{ij},X_{kl}-\imath Y_{kl}]=0,\quad1\leq i,j,k,l\leq\ell+1\,.
 \ee
Then by  \eqref{HGLsubalg}, we have
 \be
  [K^{L,R}_{ij}\,,X_{kl}-\imath Y_{kl}]\subseteq\CA\,,
 \ee
and from \eqref{HSPrelL},\eqref{HSPrelR} we deduce
 \be
  [b^L_{ij}-c^L_{ij}\,,X_{kl}-\imath Y_{kl}]\subseteq\CA\,,\quad
  [b^R_{ij}-c^R_{ij}\,,X_{kl}-\imath Y_{kl}]\subseteq\CA\,,
 \ee
which entails that $\fk\subset\hsp_{2\ell+2}$ is indeed a subalgebra
and $\CA\subset\fk$ is a Lie ideal. Thus, the spherical Lie algebra
$\fk$ is isomorphic to the semidirect sum of
$\fu_{\ell+1}\oplus\fu_{\ell+1}$ and $\CA$. $\Box$

We refer to $\fk\subset\hsp_{2\ell+2}$ as the spherical subalgebra.

\begin{de} Given spherical principal series $\hsp_{2\ell+2}$-representation
$(\pi_{\g_L,\g_R,c},\CV_{\g_L,\g_R,c})$, we define the spherical
vector $\phi\in\CV_{\g_L,\g_R,c}$ by the condition
 \be\label{SPsphvec}
  \pi_{\g_L,\g_R,c}(\fk)\cdot\phi\,=\,0\,.
 \ee
\end{de}
Note that the $\hsp_{2\ell+2}$-spherical vector \eqref{SPsphvec}
extends the $\hgl_{\ell+1}$-spherical vector defined in \eqref{WVL}
due to
 \be
  \fk\,\cap\,\hgl_{\ell+1}\,
  =\,\Span\{K^{L,R}_{ij}\,,i<j\,;\,X_{ij}-\imath Y_{ij}\,,1\leq
  i,j\leq\ell+1\}\,,
 \ee
while the consistency of \eqref{SPsphvec} is provided by Lemma
\ref{SPLEM}. Therefore,  \eqref{SPsphvec} allows the same solution
as in \eqref{LRsph}:
 \be\label{SPsph}
  \phi(x)\,=\,\varphi_{\g_L}\otimes\varphi_{\g_R}\otimes\,
  e^{-\frac{1}{2c}\,\Tr(x^\top x)}\,,
 \ee
where $\varphi_{\g_{L,R}}\in\CV_{\g_{L,R}}$ are the
$U_{\ell+1}$-spherical vectors in the $Sp_{2\ell+2}$-principal
series representations $(\pi_{\g_{L,R}} ,\CV_{\g_{L,R}})$ normalized
by $\<\varphi_{\g_{L,R}},\,\varphi_{\g_{L,R}}\>=1$.

\begin{de} Given the contragredient spherical principal series
$\hsp_{2\ell+2}$-representation
$(\pi^{\vee}_{\g_L,\g_R,c},\CV^{\vee}_{\g_L,\g_R,c})$, we define the
Whittaker covector $\psi\in\CV^{\vee}_{\g_L,\g_R,c}$ by
 \be\label{SPwhitvec}
  \pi^{\vee}_{\g_L,\g_R,c}(X_{ij})\cdot\psi\,
  =\,\imath\delta_{ij}\,\psi,\\
  \pi^{\vee}_{\g_L,\g_R,c}(a^L_{ij}+a^R_{ij})\cdot\psi\,
  =\,\pi^{\vee}_{\g_L,\g_R,c}(b^L_{ij}-c^R_{ij})\cdot\psi
  =\,\pi^{\vee}_{\g_L,\g_R,c}(b^R_{ij}-c^L_{ij})\cdot\psi=0\,.
 \ee
\end{de}
The relations \eqref{SPwhitvec} naturally extend the defining
relations \eqref{WVR} for the $\hgl_{\ell+1}$-Whittaker covector,
since the consistency of \eqref{SPwhitvec} follow by
\eqref{HSPrelL},\eqref{HSPrelR}. The equations \eqref{SPwhitvec}
allow the same solution as \eqref{WVR}. Thus  imposing $\g_L+\g_R=0$
we denote $\g_R=\g$ and keep the notation $(\pi_{\g,c},\CV_{\g,c})$
for the representation $(\pi_{-\g,\g,c},\CV_{-\g,\g,c})$, so the
equations \eqref{SPwhitvec} allow the solution \eqref{Whitvec},
 \be\label{SPwhit}
  \psi(x)\,=\,\Id_{\End(\CV_{\g})}
  \otimes\prod_{i,j=1}^{\ell+1}\delta(x_{ij}-\delta_{ij})\,,
 \ee
where $\Id_{\End(\CV_{\g})}\,\in\,\End(\CV_{\g})
  \simeq\CV_{\g}^{\vee}\otimes\CV_{\g}$
is  the identity endomorphism in the algebra $\End(\CV_{\g})$.

Thus the $\CH GL_{\ell+1}$-matrix elements $W_{\g,c}(g_1,g_2)$
defined in \eqref{GUBO},\eqref{GUBO0} can be lifted to the $\CH
Sp_{2\ell+2}$-matrix element. In particular, \eqref{GUBO0} takes the
form, for $g_1\otimes g_2\in Sp_{2\ell+2}\times Sp_{2\ell+2}$:
 \be\label{SPmatel}
  W_{\gamma,c}(g_1,g_2)\,
  =\,\<\psi,\,\pi_{\gamma_1,\g_2,c}(g_1\otimes g_2)\,\phi\>\\
  =\!\int\limits_{\Mat_{\ell+1}}\!\!\!dx
  \prod_{i,j=1}^{\ell+1}\delta(x_{ij}-\delta_{ij})\,
  \Phi_{\g}(g_1^{-1}g_2)\,e^{-\frac{1}{2c}\Tr[(g_1^{-1}xg_2)^{\top}g_1^{-1}xg_2]}\\
  =\Phi^{\ssp_{2\ell+2}}_{\gamma}(g_1^{-1}g_2)\,\,
  e^{-\frac{1}{2c}\Tr[(g_1^{-1}g_2)^{\top}g_1^{-1}g_2]}\,,
 \ee
where $\Phi^{\ssp_{2\ell+2}}_{\gamma}(g)$ is the zonal spherical
function of the spherical principal series
$\ssp_{2\ell+2}$-representation $(\pi_{\g},\CV_{\gamma})$.

\noindent {\small {\bf A.A.G.} {\sl Laboratory for Quantum Field
Theory
and Information},\\
\hphantom{xxxx} {\sl Institute for Information
Transmission Problems, RAS, 127994, Moscow, Russia};\\
\hphantom{xxxx} {\it E-mail address}: {\tt anton.a.gerasimov@gmail.com}}\\
\noindent{\small {\bf D.R.L.} {\sl Laboratory for Quantum Field
Theory
and Information},\\
\hphantom{xxxx}  {\sl Institute for Information
Transmission Problems, RAS, 127994, Moscow, Russia};\\
\hphantom{xxxx} {\it E-mail address}: {\tt lebedev.dm@gmail.com}}\\
\noindent{\small {\bf S.V.O.} {\sl
 Beijing Institute of Mathematical Sciences and Applications\,,\\
\hphantom{xxxx} Huairou District, Beijing 101408, China};\\
\hphantom{xxxx} {\it E-mail address}: {\tt oblezin@gmail.com}}

\end{document}